\documentclass{pramana}

\usepackage{graphicx,amsmath,bm}
\usepackage{url}

\begin{document}

\title{New Power Method for Solving Eigenvalue Problems}


\author{I Wayan Sudiarta\textsuperscript{1,*}\and Hadi Susanto\textsuperscript{2}}
\affilOne{\textsuperscript{1} Physics Study Program, Faculty of Mathematics and Natural Sciences, University of Mataram, Mataram, NTB, Indonesia 83125\\}
\affilTwo{\textsuperscript{2} Khalifa University, Abu Dhabi Campus,Abu Dhabi UAE, PO Box 127788} 


\twocolumn[{

\maketitle

\corres{wayan.sudiarta@unram.ac.id}

\msinfo{1 January 2024}{1 January 2024}{1 January 2024}

\begin{abstract}
We present a new power method to obtain solutions of eigenvalue problems. The method can determine not only the dominant or lowest eigenvalues but also all eigenvalues without the need for a deflation procedure. The method uses a functional of an operator (or a matrix) to select or filter an eigenvalue. The method can freely select a solution by varying a parameter associated to an estimate of the eigenvalue. The convergence of the method is highly dependent on how closely the parameter to the eigenvalues. In this paper, numerical results of the method are shown to be in excellent agreement with the analytical ones. 
\end{abstract}

\keywords{Power Iteration, Functional of Matrix, Schr\"odinger Equation}

\pacs{02.60.-x ; 03.65.Ge; 03.67.Lx; 03.65.-w;   02.70.Hm}

}]


\doinum{12.3456/s78910-011-012-3}
\artcitid{\#\#\#\#}
\volnum{123}
\year{2024}
\pgrange{23--25}
\setcounter{page}{23}
\lp{25}

\section{Introduction}
Eigenvalue problems (EVPs), which are found in many fields such as mathematics, sciences, and engineering, continue to be actively studied due to remaining open problems \cite{golub2000} and diverse applications. In physics and chemistry, the time-independent Schr\"odinger equation is an EVP where the the ground state energies and few excited states are computed \cite{levine2009, mcardle2020, adelman2021}. In modern applications, particularly for data analysis, solving EVPs is used in connection with the principal component analysis (PCA) to obtain the main patterns or dominant data \cite{demvsar2013, jolliffe2016}. 

In recent years, quantum computers are being developed and gaining acceptance to be the future of hardware for computations \cite{gyongyosi2019}.  A quantum computer uses the principles of quantum mechanics in its operation, therefore it requires quantum algorithms. By exploiting entangle properties of qubits systems, one might hope that much faster computation can be performed in comparison with the classical computers. In-line with this, new algorithms for solving eigenvalue problems on quantum computers are also being studied \cite{kandala2017, wei2019, cao2019, kyriienko, huggins2020, bauman2020, he2020, bharti2021}.

EVPs may take the form of ordinary or partial differential equations, as well as matrix equations. For computations, a differential equation can be transformed into a matrix equation by employing basis functions \cite{jackson2006} or discrete representations or grid-based approaches \cite{mathews1992}.  

Numerical methods for solving a general EVP to obtain eigenvalues and eigenvectors depend on the type, size of matrices, and the number of desired solutions. Numerous methods have been developed in the past, including the power method, Lanczos method, Arnoldi method, and QR algorithm \cite{golub2013}. Other methods can be found in various textbooks such as those by Watkins \cite{watkins2004}, Saad \cite{saad2011}, Stewart \cite{stewart2001} and Golub and Van Loan \cite{golub2013}. An overview of eigenvalue computations is provided by Golub \cite{golub2000}.  

For small-sized symmetric matrices, the QR algorithm is employed to obtain all eigen solutions. The power method or the inverse power method can be used to obtain the few dominant eigenvalues or the low-laying eigenvalues with a deflation procedure. The shifted inverse power method can select one eigenvalue in each iteration provided the chosen shift value is near the desired eigenvalue. Additionally, the power method with an exponential operator form, commonly used in imaginary time simulation, can be applied to get the few lowest eigenvalues of quantum mechanical systems. Further details can be found in \cite{sudiarta2007, sudiarta2018non, sudiarta2018}.      

For systems with a large number of variables or a large-sized matrix, methods like the shifted inverse power method, polynomial filtering, and rational filtering methods are generally preferred especially to obtain few eigenvalues \cite{saad2011, zhou2006, fang2012, li2016, xi2016}. To apply the shifted inverse power method and rational filtering method, approximately $2N^3/3$ arithmetic operations \cite{mathews1992} for solving linear equations in each step. This can be time-consuming for large matrices. Even worst requirement is when inversion of the matrix is performed before using the power method. This problem can be addressed by approximating the inverse operator by Fourier approximation \cite{childs2017, kyriienko, he2022}. However, achieving accurate approximation necessitates to add many terms, which may slow the iteration. If few interior eigenvalues are needed, instead of the  inverse power method \cite{mathews1992, golub2000}, one can use other functional or polynomial filtering method. The Chebyshev polynomials of the first kind have been widely used for approximating the $\delta$-Dirac function located near the desired eigenvalue \cite{fang2012, li2016, li2019}. However, this approach becomes ineffective when the range of eigenvalues is exceptionally large, as it requires transforming the eigenvalue problem such that all eigenvalues lie in the interval $[-1,1]$ and estimating of the lowest and highest eigenvalues.

In this paper, we introduce a novel filtering approach that enables obtaining any eigenvalues without requiring the transformation of the eigenvalue problem or estimating its lowest and highest eigenvalues, as done in Chebyshev polynomial methods. The basic principle of our filtering method is an application of a product of the matrix and the exponential of the matrix. This matrix functional  amplifies the corresponding eigenvectors related to the desired eigenvalue. Similar to the shifted inverse power method, an initial estimation of the target eigenvalue is necessary for initiating the iterative process.

The main purpose of this paper is to establish the validity of the new power method, which, to the best of our knowledge, has not been previously proposed. This paper does not intend to compare or compete with existing numerical techniques. Instead, this new approach can be considered as an additional tool for solving eigenvalue problems. Its simplicity makes it potentially useful in various applications within mathematics and the physical sciences, especially when dealing with large-sized matrices or systems.

The remaining parts of this paper is structured into four sections that explain: (i) the theoretical derivation of the new power method, (ii) numerical methods, (iii) computational results for various cases, starting with a simple matrix equation and ending with a three dimensional quantum system, and (iv) the conclusions.

\section{Theoretical Background}
The equation for eigenvalue problems is in the following form,
\begin{equation}
\hat{H}\psi = E\psi .
\label{eqn-eigen}
\end{equation} 
Notations or symbols used in Eq.\ (\ref{eqn-eigen}) are generally found for quantum systems where $\hat{H}$ is a Hamilton operator (or a matrix) and $\psi$ is a wavefunction (or a vector). The utilization of quantum mechanical notations in this paper should not lead to confusion, as they can be readily expressed in linear algebraic forms without any difficulty. To simplify notations and explanations in this paper, without loss of generality, only the operator-wavefunction form is discussed. 

We assume that the solutions of Eq.\ (\ref{eqn-eigen}) are in the form of eigenvalues $\{E_n\}$ and eigenfunctions $\{\phi_n(\mathbf{r})\}$  where $n = 0, 1, \ldots, N$. The spatial variable $\mathbf{r}$ can be in any dimensions. The number of solutions can be very large or infinite. We also assume that the eigenvalues are all reals, distinct, ordered and positives. The eigenfunctions are assumed to be orthogonal and normalized. We can index the eigenvalues according to their order,  $ 0 < E_0 < E_1 < E_2 < \ldots < E_N$. $E_0$ and $E_N$ are the lowest and the largest eigenvalues respectively. Positive eigenvalues requirement can be achieved by shifting the eigenvalues using a modified operator such as $\hat{H} + \sigma \hat{I}$ where $\hat{I}$ is the identity operator and $\sigma$ is the shift value. The requirement is important to get a stable iteration especially when an approximate operator exponentiation is used.

The power method can generally be used to determine the largest eigenvalue, in this case $E_N$. The power method begins by using an initial function, $\psi^{(0)}(r)$ and then apply the Hamilton operator to the initial function. For derivation of the power method, the initial function is expanded into the eigenfunctions given by
\begin{equation}
\psi^{(0)}(r) = \sum_{n=0}^{N} c_n \phi_n(r).
\label{eqn-initfunc}
\end{equation} 

Applying repeatedly the Hamilton operator many times (or the power operation) to the initial function, Eq.\ (\ref{eqn-initfunc}), the following function is obtained
\begin{equation}
\psi^{(k)}(r) = \hat{H}^k \psi^{(0)}(r) = \sum_{n=0}^{N} c_n \hat{H}^k \phi_n(r) = \sum_{n=0}^{N} c_n E_n^k \phi_n(r).
\label{eqn-power1}
\end{equation}
After a large number of iterations, $K$, it is shown in Eq.\ (\ref{eqn-power1}) that the largest value or dominant term is for the largest eigenvalue term. The wavefunction can then be approximated by  
\begin{equation}
\psi^{(K)}(r) \approx [c_N E_N^K] \phi_N(r).
\label{eqn-power2}
\end{equation}
The factor in the square bracket can be eliminated by a normalization procedure. For numerical purposes, in order to avoid overflow in calculations,  a scaling procedure or a normalization of the wavefunction is performed in every iteration. The largest eigenvalue is obtained provided that the initial function contains the largest eigenvalue wavefunction or its coefficient is non-zero, $c_N \neq 0$.

For every iteration step, an estimate of the eigenvalue can be determined by using
\begin{equation}
E^{(k)} = \frac{\int \psi^{(k)*}\hat{H} \psi^{(k)} dV}{\int \psi^{(k)*} \psi^{(k)} dV},
\end{equation}
where $dV$ is the volume element and the interval of integration is set appropriately according to the spatial domain of the problem.

The power method can also be applied to get the lowest eigenvalue. One can use the inverse power method that applies an inverse of the operator or requires solving a system of linear equations. Another method is by using a function of the Hamilton operator.  The function is chosen such that it amplifies the lowest eigenvalue. For applications in quantum systems, applying the finite difference time domain (FDTD) method \cite{sudiarta2007} and Monte Carlo (MC) method \cite{foulkes2001quantum, toulouse2016introduction}, the exponential function is often used. The eigen equation (Eq.\ \ref{eqn-eigen})  or the Schr\"odinger equation can be transformed to a diffusion equation which then be solved by an evolution method \cite{sudiarta2007}.  The exponential of the Hamilton operator is given by
\begin{equation}
\hat{G} = e^{-\tau \hat{H}}
\label{eqn-G}
\end{equation}
where $\tau $ is a constant related to the interval of simulation time. Similarly, by repeatedly applying the operator $\hat{G}$ many times to the initial function, a new function is obtained, i.e.,
\begin{eqnarray}
\psi^{(k)}(r) &=& \hat{G}^k \psi^{(0)}(r) = \sum_{n=0}^{N} c_n \left[e^{-\tau \hat{H}}\right]^k \phi_n(r) \nonumber \\
&=& \sum_{n=0}^{N} c_n e^{- k\tau E_n} \phi_n(r).
\label{eqn-power3}
\end{eqnarray}
For a large value of $k$ or $K$, it is shown in Eq.\ (\ref{eqn-power3}) that the dominant term is the lowest eigenvalue term.  The wavefunction approaches the lowest eigenvalue wavefunction
\begin{equation}
 \psi^{(K)}(r) \approx [c_0 e^{-K \tau E_0}] \phi_0(r).
\label{eqn-power4}
\end{equation}

For determining higher eigenvalue (or excited state) wavefunctions, the same iteration procedure is performed but the wavefunction $\psi$ must be orthogonal to the lowest eigenvalue wavefunction. This procedure is known as explicit deflation method \cite{golub2000}.

In this paper, we propose a new power method to obtain any eigenvalue without the deflation procedure. The new power method utilizes a product of the operators $\hat{H}$ and $\hat{G}$. The new power method uses the following  operator:
\begin{equation}
\hat{F} = \hat{H}^\alpha e^{-\tau \hat{H}}.
\end{equation}

In the same process as the power method, applying $\hat{F}$ operator many times to the initial function, the wavefunction after $k$th iteration is given by
\begin{eqnarray}
\psi^{(k)}(r) &=& \hat{F}^k \psi^{(0)}(r) = \sum_{n=0}^{N} c_n \left[\hat{H}^\alpha e^{-\tau \hat{H}}\right]^k \phi_n(r) \nonumber \\ &=& \sum_{n=0}^{N} c_n \left[f(E_n)\right]^k \phi_n(r)
\label{eqn-power5}
\end{eqnarray}
where the function $f(E)$ is given by
\begin{equation}
f(E) =  E^{\alpha} e^{-\tau E}. 
\label{eqn-fe}
\end{equation}

\begin{figure}[htp]
\begin{center}
\includegraphics[width=\columnwidth]{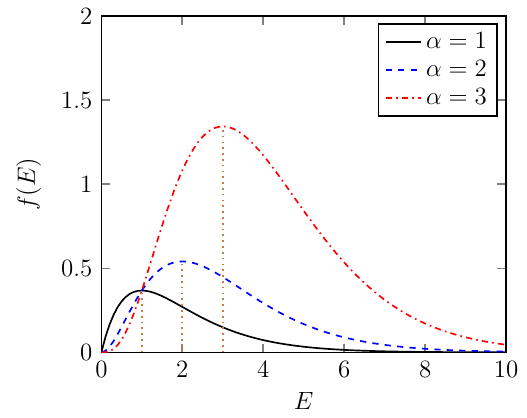}
\end{center}
\caption{$f(E)$ in Eq.\ (\ref{eqn-fe}) as a function of $E$ with $\tau = 1$ for three values of $\alpha$. The dotted brown lines are used to indicate the positions of the maximum values.}
\label{fig-fe}
\end{figure}

Examples of $f(E)$ curves as a function of $E$ are shown in Fig.\ \ref{fig-fe}. It is noted that $f(E)$ in the interval $[0,\infty]$ has one peak located at $E_{p} = \alpha/\tau $. The value of $\tau$ can be chosen appropriately such that $E_{p}$ is closed to a desired eigenvalue. Substituting $\tau = \alpha/E_{p}$ into $f(E)$ (Eq.\ (\ref{eqn-fe}), it is found that $f(E) = [E \exp(-E/E_{p})]^\alpha$. This result shows that the value of $\alpha$ does not produce a significant improvement on the curve $f(E)$ since it is the same as the $\alpha$th power of $f(E)$ (with $\alpha = 1$).  Therefore, in the rest of the paper, $\alpha = 1$ is used.

By choosing an appropriate value of $E_{p}$, we can select one solution with an eigenvalue that is near to $E_{p}$. This means that after large number iterations, $k = K$,
\begin{equation}
\psi^{(K)}(r) \approx [c_n f(E_n)^K ] \phi_n(r).
\label{eqn-power6}
\end{equation}
where the index $n$ is for the maximum value of $f(E_n)$.
 
This iteration procedure is useful to obtain an eigensolution when its approximate value or $E_p$ is known.  The value of $E_{p}$ can be also varied in order to obtain all eigensolutions.
 
For a chosen value of $E_p$, three eigenvalues near $E_p$ are assumed to be $E_a$, $E_b$ and $E_c$. These eigenvalues are ordered such that $f(E_a) > f(E_b) > f(E_c)$. Following Golub and Van Loan \cite{golub2013}, the convergence of the method is determined by the ratio, $R(E_p)$, given by
 \begin{equation}
 R(E_p) =  \frac{f(E_b)}{f(E_{a})}  = \frac{E_b}{E_{a}} e^{-(E_b-E_{a})/E_{p}} 
 \end{equation}

The error of the computed eigenvalue after $k$th iteration is given by
\begin{equation}
\text{Error}^{(k)} = |E_n - E^{(k)}| \approx  \beta |R(E_p)|^k
\end{equation}
where $\beta$ is a proportionality constant.

If the iteration reaches a convergent eigenvalue after the error below a certain tolerance value, $tol$, the number iteration needed is
\begin{equation}
k  = \frac{C}{\log(|R(E_p)|)}
\label{eqn-kval}
\end{equation}
where a contant $ C = \log(tol/\beta)$.
 
In Eq.\ (\ref{eqn-kval}), it is noted that the number of iterations becomes very large (or infinite) or the convergence is very slow when $R(E_p) \approx 1$. This happens when  $E_p$ equals to the turning point value $E_{tp}$ given by
\begin{equation}
E_{tp} = \frac{E_b - E_a}{\ln(E_b/E_a)}.
\label{eqn-tp}
\end{equation}

\section{Numerical Methods}
\subsection{Matrix Form}
For numerical computations, the short-time approximation of exponential of an operator (a matrix) is used \cite{trotter1959, suzuki1976, lloyd1996, burkard2022} as follows:  
\begin{equation}
e^{-\tau \hat{H}} = \left[e^{-\tau \hat{H}/M}\right]^M  
\end{equation}

Applying the first order Taylor approximation $e^{-\tau\hat{H}/M} \approx (1 - \frac{\tau}{M}\hat{H}) $, it is obtained
\begin{equation}
e^{-\tau \hat{H}}  \approx \left[ 1 - \frac{\tau}{M}\hat{H} \right]^M. 
\label{eqn-exp1}
\end{equation}

Therefore,  the operator used in the new power method is 
\begin{equation}
\hat{F} \approx \hat{H}\left[ 1 - \frac{1}{M E_p}\hat{H} \right]^M 
\label{eqn-shorttime}
\end{equation}

It is noted in Eq.\ (\ref{eqn-shorttime}) that besides the initial wavefunction, $\psi^{(0)}$, two other parameters that affect the performance of the iteration method are the estimate value of eigenvalue, $E_p$, and the parameter $M$. The resulting operator $\hat{F}$ in Eq.\ (\ref{eqn-shorttime}) is a polynomial of the Hamilton operator and without other additional constant operators. The resulting operator $\hat{F}$ shares the same eigenfunctions as for the operator $\hat{H}$. Therefore the power iteration with $\hat{F}$ always converges to one of the eigenfunctions of $\hat{H}$. 

\subsection{Differential Form}
For the Hamilton operator in the form of differential operators, we represent functions by discrete values on grid points, as an example in one dimensional space $\psi_i \equiv \psi(i \Delta x)$ with the spatial grid spacing $\Delta x$.  We then use finite difference approximation formulas for derivatives. For the second order derivative operator in one dimension, $\hat{D}_2 = d^2/dx^2$ is approximated \cite{mathews1992} by 
\begin{eqnarray}
\hat{D}_2 \psi(x)|_i &\approx & \frac{\psi(x+\Delta x) - 2\psi(x) + \psi(x-\Delta x)}{(\Delta x)^2} \nonumber \\  &=& \frac{\psi_{i+1} - 2\psi_i + \psi_{i-1}}{(\Delta x)^2} .
\end{eqnarray} 

For applications to the Schr\"odinger equation, the exponentiation of the operator can be performed by solving the diffusion equation \cite{sudiarta2007, sudiarta2018non} as follows.
\begin{equation}
\frac{\partial \psi(r,\tau)}{\partial \tau} = -\hat{H}\psi(r,\tau)
\end{equation}

Using the forward scheme of finite difference formula for first order derivative in time, we get similar expression as for the matrix case (see Eq.\ (\ref{eqn-exp1}), that is
\begin{equation}
\psi(r,\tau+\Delta \tau) = \psi(r, \tau)  - \Delta \tau \hat{H}\psi(r,\tau)
\label{eqn-expdif}
\end{equation}
where $\Delta \tau$ is the temporal interval. Substituting $\Delta \tau = \tau/M$ into Eq.\ (\ref{eqn-expdif}) gives a similar expression to Eq.\ (\ref{eqn-exp1}).

\section{Results and Discussions}

In this section, we give applications of the new power method to solve a simple eigenmatrix equation and the Schr\"odinger equations for a particle in (a) a one dimensional box, (b) a ring, (c) a one dimensional harmonic potential, and (d) a three dimensional cubic box. Computer codes for computing results in this paper are available in \cite{sudiarta2022}.

\subsection{A Simple Matrix}

For eigenvalues of a simple matrix, following an example in \cite{wikipedia}, we have used a $3\times 3 $ matrix to test our numerical procedure. The matrix is given by
 \begin{equation}
\hat{H} = 
\begin{bmatrix}
2 & 0 & 0\\
0 & 2 & 1 \\
0 & 1 & 2 
\end{bmatrix}
\label{eqn-splmatrix}
\end{equation}
where its eigenvalues are $1$, $2$ and $3$ and corresponding eigenvectors are $\phi_1 = (1/\sqrt{2})[0,1,-1 ]^T$, $\phi_2 = [1, 0, 0]^T$, $\phi_3 =(1/\sqrt{2}) [0, 1, 1]^T$.

\begin{figure}[htb]
\begin{center}
\includegraphics[width=\columnwidth]{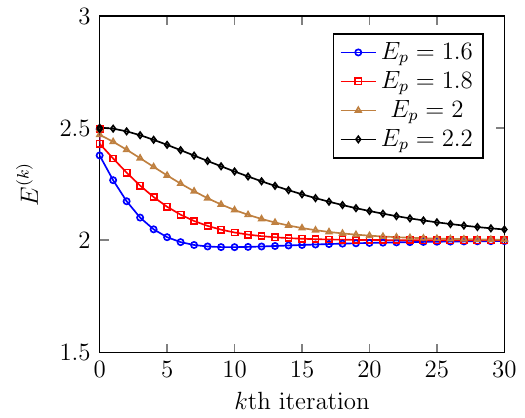}
\end{center}
\caption{Numerical values of $E$ for the matrix in Eq.\ (\ref{eqn-splmatrix}) at $k$th iteration with the iteration parameters $M = 100$ and $\psi^{(0)} = [0.7, 0.8, 0.4 ]^T$ for four values of $E_p$. }
\label{fig-fe1}
\end{figure}

\begin{figure}[htb]
\begin{center}
\includegraphics[width=\columnwidth]{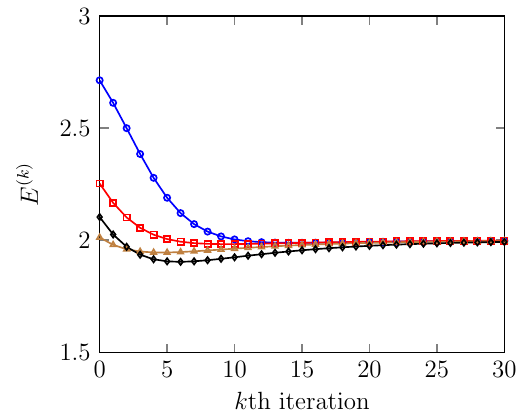}
\end{center}
\caption{Numerical values of $E$ for the matrix in Eq.\ (\ref{eqn-splmatrix}) at $k$th iteration with the parameters $M= 100$ and $E_p = 1.6$ for four different random initial vectors. }
\label{fig-fe2}
\end{figure}

\begin{figure}[htb]
\begin{center}
\includegraphics[width=\columnwidth]{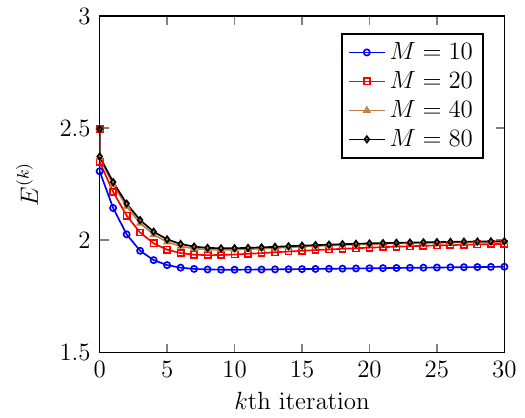}
\end{center}
\caption{Numerical values of $E$ for the matrix in Eq.\ (\ref{eqn-splmatrix}) at $k$th iteration with the parameters $E_p= 1.6$ and $\psi^{(0)} = [0.7, 0.8, 0.4]^T$  for four values of $M$. }
\label{fig-fe3}
\end{figure}

\begin{figure}[htb]
\begin{center}
\includegraphics[width=\columnwidth]{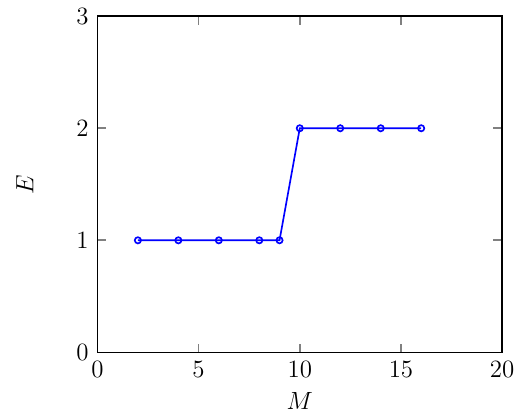}
\end{center}
\caption{Eigenvalues of the matrix in Eq.\ (\ref{eqn-splmatrix}) computed using the iteration parameters $E_p = 1.6$ and $\psi^{(0)} = [0.7, 0.8, 0.4]^T$ for various values of $M$. }
\label{fig-fe4}
\end{figure}

Before applying the new power iteration procedure, three parameters that need to be specified are the initial vector $\psi^{(0)}$, the location of the peak $E_p$ and $M$. We first study the effects of these three parameters using the simple matrix in Eq.\ (\ref{eqn-splmatrix}). For this case, we show only results for determining the second eigenvalue. Other eigenvalues are found in similar manner. Numerical results of the estimates of eigenvalues after $k$th iteration are shown in Fig.\ \ref{fig-fe1} - \ref{fig-fe3}.  Almost all parameters used in the calculations give convergence to the eigenvalue of $2$ as expected since the values of $E_p$ near $2$ are used. There is one set of numerical results for the parameter $M = 10$ in Fig.\ \ref{fig-fe3} to show a slow convergence. Further investigation of the effect due to the parameter $M$ on the resulting eigenvalue is shown in Fig.\ \ref{fig-fe4}. The numerical results show that the parameter $M$ does affect the resulting eigenvalue. Figure \ref{fig-fe4} indicates that for $M < 10$, the iteration converges to the eigenvalue of $1$ and  for $M \geq 10$ the iteration converges to the eigenvalue of $2$. This shows for this case that the short-time approximation of the matrix exponentiation is less accurate for the parameter $M < 10$. However, all numerical results give one of the eigenvalues of the simple matrix.       

\subsection{Particle in a 1D Box}

\begin{figure}[htb]
\begin{center}
\includegraphics[width=\columnwidth]{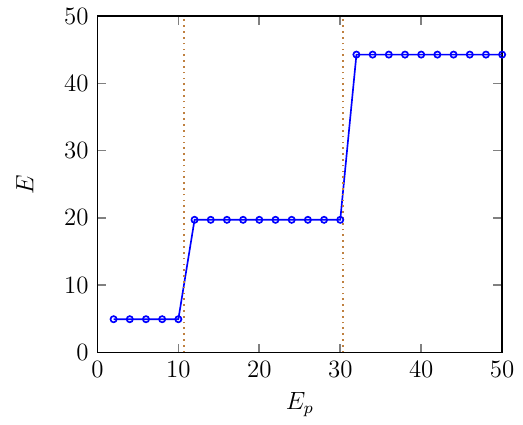}
\end{center}
\caption{Eigenvalues of a particle in a one dimensional box with a length $L = 1$ computed by the power iteration method for various values of $E_p$. The dotted lines are used to indicate the turning points $E_{tp}$ in Eq.\ (\ref{eqn-tp}).}
\label{fig-fdtd1}
\end{figure}

For the following cases in solving the Schr\"odinger equation, for simplicity in presentation, we have used natural atomic  units. For a particle in a one dimensional box with a length $L = 1$ at interval $[0,1]$,  the Hamilton operator is given by 
\begin{equation}
\hat{H} =  -\frac{1}{2}\frac{d^2}{dx^2}
\end{equation}
The eigenvalues of the system are $E_n = n^2\pi^2/2$ and the corresponding eigenfunctions are $\phi_n(x) = \sqrt{2}\sin(n\pi x)$. In order to ensure the initial function contains the ground state (or the lowest eigenvalue), a random initial wavefunction is used. The  computational parameters used in this case are $\Delta x = 0.02$, $\Delta \tau = (\Delta x)^2/10$ and $tol = 10^{-8}$. The temporal interval is chosen such that $\Delta \tau \leq (\Delta x)^2$ in order to get stable solutions of the diffusion equation. Three eigenvalues of this system are given in Table \ref{tbl-fdtd1} and Fig.\ \ref{fig-fdtd1} for various values of $E_p$. Fig.\ \ref{fig-fdtd1} shows that the resulting convergent eigenvalues are observed to be dependent on the value of $E_p$. The turning point values $E_{tp}$ are also found to be in good agreement with the analytical expression (Eq.\ (\ref{eqn-tp})).

\begin{table}[htb]
\begin{center}
\caption{Exact and numerical eigenvalues for a particle in a one dimensional box of length $L=1$.}
\begin{tabular}{crrr}
\hline
$n$ &  Exact $E_n$ & Numerical $E_n$ & $|\Delta E|/E (\%)$\\
\hline 
1	&	4.934802	&	4.933179	&	0.03	\\  
2	&	19.739208	&	19.713247	&	0.13	 \\
3	&	44.413219	&	44.281873	&	0.30	 \\
4	&	78.956835	&	78.542094	&	0.53	 \\
5	&	123.370055	&	122.358708	&	0.82	 \\
\hline
\end{tabular}
\label{tbl-fdtd1}
\end{center}
\end{table}

\begin{figure}[htb]
\begin{center}
\includegraphics[width=\columnwidth]{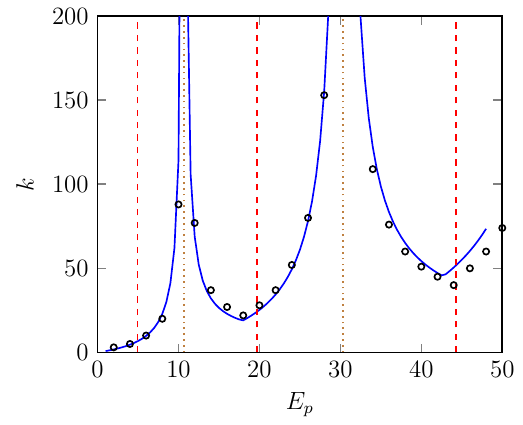}
\end{center}
\caption{Number of iterations needed to reach convergent eigenvalues of a particle in a one dimensional box with a length $L = 1$ for various values of $E_p$:  (a) numerical results (circles) (b) fitting with Eq.\ \ref{eqn-kval} using the fitting parameter $C$.  The dotted lines are used to indicate the turning points $E_{tp}$. The dashed lines indicate the positions of eigenvalues.}  
\label{fig-fdtd2}
\end{figure}

\subsection{Particle in a Ring}

As a further test to the new power method, we consider a problem of a particle moving in a ring with the circumference length $L$.  This problem can be viewed as a particle in one dimension with a periodic boundary condition (PBC). The boundaries in the previous example for one-dimensional box is replaced by the PBC which is $\psi(0,t) = \psi(L,t)$. The eigenvalues of a particle in the ring are $2n^2/L^2$ with $n = 0, \pm 1, \pm 2, \ldots$. It is noted here that there is an eigenvalue of zero for $n = 0$ and degenerate eigenvalues (with more than one solution for every eigenvalue) for $n \neq 0$. The computational parameters used here are the same as in the previous one-dimensional box. Numerical results are shown in Table \ref{tab-ring}. It is shown that the zero eigenvalue cannot be reached in the iteration due to the function $f(E)$ is zero at $E = 0$ (see also Fig.\ \ref{fig-fe}). The exclusion of the zero eigenvalue can be an advantage or a disadvantage of the method depending on whether the zero eigenvalue is important in its applications. In the case of zero eigenvalue, the operator $\hat{G}$ in Eq.\ (\ref{eqn-G}) can be used instead. The power method is based on the dominant eigenvalue and it cannot separate two or more solutions with the same eigenvalue. In that case, degenerate solutions can be separated by performing perturbation in the Hamilton operator or introducing symmetry requirements.

\begin{table}[htb]
\begin{center}
\caption{Eigenvalues of a particle in a ring with the circumference length of $L = 1$.}
\label{tab-ring}
\begin{tabular}{crrr}
\hline
n	&	Exact E 	&	Numerical E	&	$|\Delta E|/E (\%)$	\\
\hline
0   &    0     & NA & NA \\
1	&	 19.739208	&	 19.732716	&	0.03	\\
2	&	 78.956835	&	 78.852987	&	0.13	\\
3	&	177.652879	&	177.127493	&	0.30	\\
4	&	315.827340	&	314.168389	&	0.52	\\
\hline
\end{tabular}
\end{center}
\end{table}

\subsection{Particle in a Harmonic Oscillator}

Another validation test of the new power method is by using the Hamilton operator for a particle in a harmonic potential given by 
\begin{equation}
\hat{H} =  -\frac{1}{2}\frac{d^2}{dx^2} + V(x)
\end{equation}
where $V(x) = \frac{1}{2} x^2$. The eigenvalues are $ (n+\frac{1}{2})$ and the wavefunctions are $\phi_n(x) = \pi^{-\frac{1}{4}}[2^n n!]^{-\frac{1}{2}} h_n(x) e^{-\frac{1}{2}x^2}$ with the Hermite polynomial $h_n(x)$. 

The present of the potential $V(x)$ in the operator requires a slight modification in the numerical procedure for the evolution operator $e^{-\tau \hat{H}}$. Following Sudiarta and Geldart \cite{sudiarta2007}, to obtain accurate results and a stable iteration, the evolution of the wavefunction is obtained by
\begin{eqnarray}
&\psi(x, \tau+\Delta \tau)& \approx   a(x)\psi(x,\tau)  - \frac{b(x)\Delta \tau}{2(\Delta x)^2}\times  \\&&\left[\psi(x-\Delta x,\tau) - 2\psi(x,\tau) + \psi(x+\Delta x,\tau)\right] \nonumber
\label{eqn-df2}
\end{eqnarray} 
with $a(x) = [1 - V(x)\Delta \tau/2]/[1 + V(x)\Delta \tau/2]$ and $b(x) = 1/[1 + V(x)\Delta \tau/2]$.

\begin{figure}[htb]
\begin{center}
\includegraphics[width=\columnwidth]{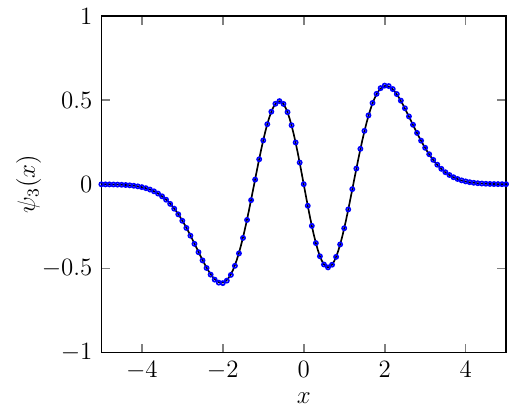}
\end{center}
\caption{numerical (circles) and exact results (line) for the third excited state wavefunction of a particle in a harmonic potential.}
\label{fig-harmonic3}
\end{figure}

The stability requirement for Eq.\ (\ref{eqn-df2}) is $b(x)\Delta \tau <(\Delta x)^2$. Assuming that the potential used is always positive, $V(x) > 0$, the coefficient $b(x)$ is always less than one, $b(x) < 1.0$. Therefore, the same stability requirement $\Delta \tau <(\Delta x)^2$ as in the previous example can also be used in this case. This is one of advantages of using Eq.\ (\ref{eqn-df2}) compared with Eq.\ (\ref{eqn-expdif}). The computational parameters used are $\Delta x = 0.1$ and $\Delta \tau = (\Delta x)^2/10$. The interval of computation is $[-10,10]$. Numerical results for the eigenvalues and the third excited state wavefunction $\phi_3(x)$ of a particle in a harmonic potential are given in Table \ref{tab-harmonic} and Fig.\ \ref{fig-harmonic3} respectively. The numerical results agree well with the exact results. The numerical errors of the eigenvalues are less than $0.5\%$. The errors are mainly due to the approximation of the derivative by the finite difference method and not due to the new power method.

\begin{table}[htb]
\begin{center}
\caption{Eigenvalues of a particle in a harmonic oscillator potential.}
\label{tab-harmonic}
\begin{tabular}{crrr}
\hline
$n$	&	Exact E 	&	Numerical E	&	$|\Delta E|/E (\%)$	\\
\hline
0	& 0.5	&	0.499687	&	0.06 \\
1	& 1.5	&	1.498437	&	0.10 \\
2	& 2.5	&	2.495937	&	0.16 \\
3 &	3.5	&	3.492195	&	0.22\\
4 &	4.5	&	4.487217	&	0.28\\
5 &	5.5	&	5.480985	&	0.35\\
6 &	6.5	&	6.473401	&	0.41\\
\hline
\end{tabular}
\end{center}
\end{table}

\subsection{A Particle in 3D Cubic Box}
In this section we apply the new power method for a three dimensional operator. Similarly to the one-dimensional case, the differential operator is approximated by the finite difference formula. For a particle in a three dimensional box with the side length $L = 1$ (i.e., with the interval $x,y,z \in [0,1]$), the Hamilton operator for this case is given by 
\begin{equation}
\hat{H} =  -\frac{1}{2}\left[\frac{\partial^2}{\partial x^2} + \frac{\partial ^2}{\partial y^2} + \frac{\partial ^2}{\partial z^2} \right]
\end{equation}
The eigenvalues are $E_{n_x n_y n_z} = (\pi^2/2) (n_x^2 + n_y^2 + n_z^2)$ and the eigenfunctions are  
\begin{equation}
\phi_{n_x,n_y,n_z}(x, y,z) = \sqrt{8}\sin(n_x\pi x)\sin(n_y\pi y)\sin(n_z\pi z).
\end{equation}

The computational parameters used are $\Delta x = 0.05$ and $\Delta \tau = (\Delta x)^2/10$. 

\begin{figure}[htb]
\begin{center}
\includegraphics[width=\columnwidth]{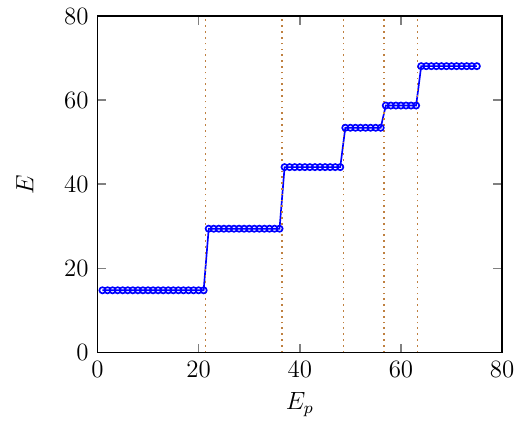}
\end{center}
\caption{Eigenvalues of a particle in a cubic box computed by the power iteration method for various values of $E_p$. The dotted lines are used to indicate the positions of $E_{tp}$ computed by Eq.\ (\ref{eqn-tp}).}  
\label{fig-cubic}
\end{figure}

Numerical results of eigenvalues for the particle in the cubic box computed by the new power method are found in Table \ref{tab-cubic}. Isosurfaces of the wavefunction $\phi_{222}(x,y,z)$ computed numerically are given in Fig.\ \ref{fig-cubic}. The numerical results are in good agreement with the exact results with the errors less than $2\%$. This corresponds to the small number of grid points per dimension (about 20) used in the computations due to a limitation in our computer memory. The turning point values in Fig.\ \ref{fig-cubic} are all found to be consistent with Eq.\ (\ref{eqn-tp}).  Using larger number of grid points or smaller spatial spacing can certainly improve the accuracy of the finite difference approximation.

\begin{table}[htb]
\begin{center}
\caption{Eigenvalues of a particle in a cubic box with side length of 1.}
\label{tab-cubic}
\begin{tabular}{crrr}
\hline
$n_x,n_y,n_z$	&	Exact E 	&	Numerical E	&	$|\Delta E|/E (\%)$	\\
\hline
1,1,1	&	14.804407	&	14.773991	&	0.21	\\
1,1,2	&	29.608813	&	29.426721	&	0.61	\\
1,2,2	&	44.413220	&	44.079451	&	0.75	\\
1,1,3	&	54.282824	&	53.446710	&	1.54	\\
2,2,2	&	59.217626	&	58.732179	&	0.82	\\
1,2,3	&	69.087231	&	68.099436	&	1.43	\\
\hline
\end{tabular}
\end{center}
\end{table}

\begin{figure}[htb]
\begin{center}
\includegraphics[width=\columnwidth]{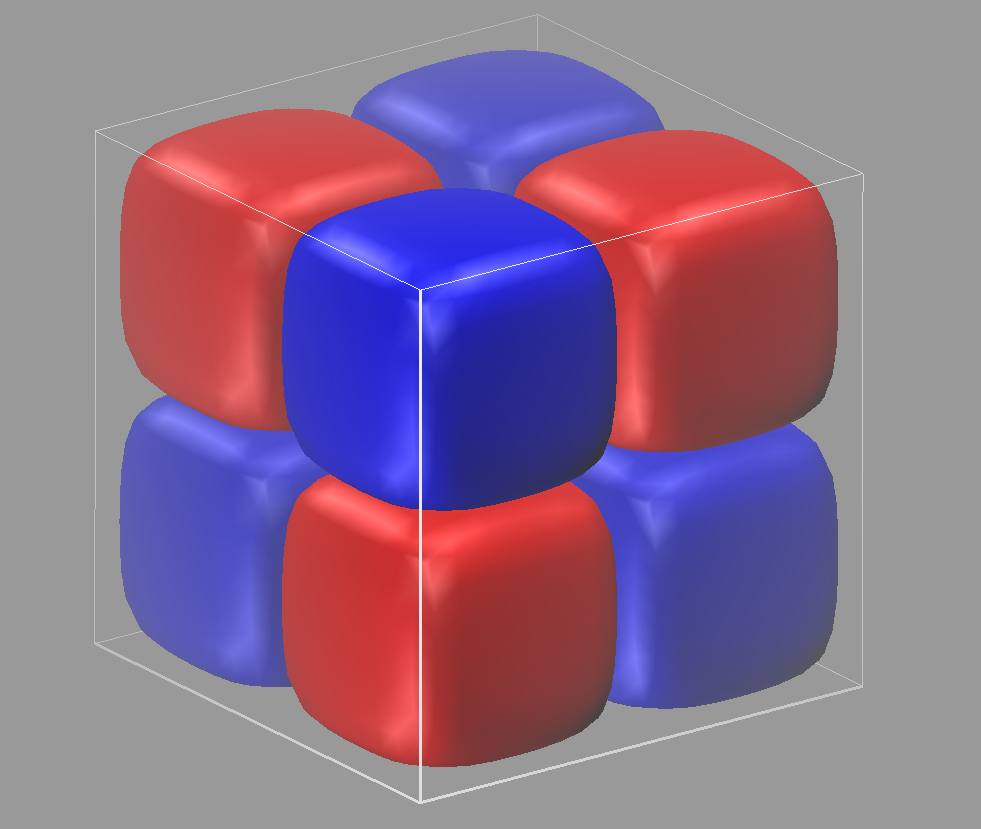}
\end{center}
\caption{Isosurfaces of the wavefunction $\phi_{222}(x,y,z)$ for a particle in a cubic box with the side length of 1. The values of isosurfaces are $0.1$ (red) and $-0.1$ (blue).}
\label{fig-isosurface3}
\end{figure}

\section{Conclusions}
We have presented a new power method with a functional of an operator to obtain solutions of eigenvalue problems. The numerical results of validation examples for various operators have shown that the new power method  always produces a convergent eigenvalue closed to a parameter (or an approximate eigenvalue) of the functional used. The convergence of the method is highly dependent on the closeness of the parameter to one of the eigenvalues of the operator. The analytical expressions for the convergence and the number of iterations needed to get a convergent eigenvalue have also been derived. The numerical results for the number of iterations have agreed well with the analytical results. The new power method can be used to select freely one solution of an eigenvalue problem by adjusting the parameter of the functional. Therefore, the new power method extends the application of power methods to obtain not only the dominant eigenvalue solution but also all eigensolutions.   

\section{Acknowledgements}
IWS is partially supported by University of Mataram through a PNBP Internal Research Grant. HS is supported by Khalifa University through a Faculty Start-Up Grant (No.\ 8474000351/FSU-2021-011) and a Competitive Internal Research Awards Grant (No. 8474000413/CIRA-2021-065).

\bibliography{references-v2}

\end{document}